\documentclass[11pt]{article}
\usepackage{latexsym}

\usepackage[tbtags]{amsmath}
\usepackage{epsfig,amstext,amssymb,amsthm,latexsym}
\usepackage{amssymb,color}

\UseRawInputEncoding

\definecolor{c20}{rgb}{0.,0.7,0.}
\definecolor{c30}{rgb}{0.,0.,1.}
\definecolor{c40}{rgb}{1,0.1,0.7}
\definecolor{c50}{rgb}{1,0,0}

\date{}
\newtheorem{theorem}{Theorem}[section]

\newtheorem{proposition}{Proposition}[section]

\newtheorem{remark}{Remark}[section]

\makeatletter 
\@addtoreset{equation}{section}
\makeatother 

\begin{document}
\title{On Rank Correlation Coefficients}

\author{ Alexei Stepanov \thanks{\noindent  Education and Research Cluster ``Institute of High Technology",\  Immanuel Kant  Baltic Federal University, A.Nevskogo 14, Kaliningrad, 236041 Russia,  email: alexeistep45@mail.ru}}

\maketitle
\begin{abstract} 
In the present paper, we propose  a new rank correlation coefficient $r_n$, which is a sample analogue of the theoretical correlation coefficient $r$, which, in turn, was proposed in the recent work of Stepanov (2025b). We discuss the properties of $r_n$ and compare $r_n$ with  known rank Spearman $\rho_{S,n}$, Kendall $\tau_n$  and   sample Pearson $\rho_n$ correlation coefficients. Simulation experiments show that when the relationship between $X$ and $Y$ is not close to linear, $r_n$ performs better than other correlation coefficients.  We also find analytically the values  of $Var(\tau_n)$ and $Var(r_n)$. This allows to estimate theoretically the asymptotic performance of $\tau_n$ and $r_n$. 
\end{abstract}



\section{Introduction} Let  $(X,Y), (X_1,Y_1),\ldots,(X_n,Y_n)$  be independent and identically distributed  random vectors with absolutely continuous bivariate distribution function $F(x,y)=P(X\leq x, Y\leq y)$,  bivariate survival function $\bar{F}(x,y)=P(X>x,Y>y)$, density function  $f(x,y)$ and marginal distribution and density functions    $H(x), G(y), h(x)$ and $g(y)$. There are three basic  measures of dependence rate: the Pearson $\rho$, the Spearman $\rho_S$ and the Kendall $\tau$ correlation coefficients. 

The rate of dependence between the variables $X$ and $Y$ is most often measured by the Pearson correlation coefficient
$$
\rho=\frac{E(X-EX)(Y-EY)}{\sigma_X\sigma_Y}.
$$
The random variable
$$
\rho_n =\frac{\sum\limits_{i=1}^n (X_i-\bar{X})(Y_i-\bar{Y})} {\sqrt{\sum\limits_{i=1}^n (X_i-\bar{X})^2 \sum\limits_{i=1}^n (Y_i-\bar{Y})^2}}
$$
represents the sample Pearson correlation coefficient, for which it is known that $\rho_n\stackrel{p}{\rightarrow}\rho$; see, for example,  Fisher (1921). 

Let  $X_{(1)}\leq\ldots\leq X_{(n)}$ be the order statistics obtained from the sample   $X_1,\ldots,X_n$. For these order statistics, let us define their concomitants   $Y_{[1]},\ldots,Y_{[n]}$. Let $X_i=X_{(j)}$, then  $Y_{[j]}=Y_i$ be the concomitant of the order statistic  $X_{(j)}$.  The concept of concomitants is   discussed in the papers of   David and Galambos (1974), Bhattacharya (1974), Egorov and Nevzorov (1984),  David and Nagaraja (2003), Balakrishnan and Lai (2009), Bairamov and Stepanov (2010),   Balakrishnan and Stepanov (2015). See also the references therein. 

Let  $I_{ji}=I(Y_{[j]}\leq Y_{[i]}),\ i\not= j\in \{1,2,\ldots,n\}$, where $I$ is the indicator function. The random variable  
$$
\tau _n=\frac{4\sum_{i=2}^n\sum_{j=1}^{i-1}I_{ji}}{n(n-1)}-1
$$
is known as the rank Kendall correlation coefficient. The theoretical analogue of $\tau_n$ is presented by
$$
\tau=4E[F(X,Y)]-1=4\int_{\mathbb{R}^2}F(x,y)f(x,y)dxdy-1.
$$
It is known, see, for example, Balakrishnan and Lai (2009) or Stepanov (2025a), that $\tau=E\tau_n$.  The random variable  
$$
\rho_{n,S} =1-\frac{6\sum_{i=1}^n\left(\sum_{j=1}^{n}I_{ji}-i\right)^2}{n^3-n}
$$ 
is known as the rank Spearman correlation coefficient. The theoretical analogue of $\rho_{n,S}$ is presented by
$$
\rho_S =12E[H(X)G(Y)]-3=12\int_{\mathbb{R}^2}H(x)G(y)f(x,y)dxdy-3.
$$ 
New theoretical and rank correlation coefficients were introduced and analyzed in Stepanov (2025b). The  theoretical correlation coefficient  $r$ was defined by 
\begin{eqnarray*}
r&=&6E[F(X,Y)-H(X)G(Y)]\\
&=& 6\int_{\mathbb{R}^2}[F(x,y)-H(x)G(y)]f(x,y)dxdy
\end{eqnarray*}
and the rank correlation coefficient $\tilde{r}_n$  was defined by $\tilde{r}_n= \frac{3\tau_n-\rho_{n,S}}{2}$. The rationale for introducing $r$ is as follows. Let us consider the dependence function $F(x,y)-H(x)G(y)$, which gives us the dependence value in every point $(x,y)$ of the distribution support. We adjust the definition of the dependence function  and redefine it as $6(F(x,y)-H(x)G(y))$. We can view any correlation coefficient     as  some average value of  the dependence rate, or of the dependence function. Observe that the dependence function does not take into account the likelihood of obtaining points $(x,y)$. In the above definition of $r$, we  take this likelihood into account. Thus,  $r$ represents  the ``ideal" average of the dependence function. It was shown in Stepanov (2025b) that
$$
r= \frac{3\tau-\rho_{S}}{2}
$$
and
$$
E\tilde{r}_n=\frac{3n\tau-(n-2)\rho_S}{2(n+1)}\rightarrow r\quad  (n\rightarrow \infty).
$$
Some other properties of $r$ and $\tilde{r}_n$ were discussed there. It should be noted that  the form of $\tilde{r}_n$ is complex. It was suggested in Stepanov (2025b) that  a simpler  rank correlation coefficient can be offered instead of $\tilde{r}_n$. This simple rank correlation coefficient $r_n$  is proposed in our present work.

The  correlation coefficients $\rho_n,\ \rho_{n,S},\ \tau_n$ and $\tilde{r}_n$ were compared in Stepanov (2025b). It  was revealed  that the behavior of $\rho_n$ can be very different from the behavior of the rank correlation coefficients $\rho_{n,S},\ \tau_n,\ \tilde{r}_n$, which, in turn, behave in a similar way with each other. This  follows from the definition of $\rho_n$ and $\rho_{n,S},\ \tau_n\ \tilde{r}_n$. The coefficient $\rho_n$ measures the association between variables $X$ and $Y$, whereas  the coefficients $\rho_{n,S},\ \tau_n$ and $\tilde{r}_n$ measure the association between their ranks. A question was raised:  which  correlation coefficient measures the dependence rate best?  Based on the simulation analysis of the variances of correlation coefficients, it was concluded in Stepanov (2025b) that  $\rho_n$ is preferable when the association between $X$ and $Y$ is "close``\ to $\pm 1$ and  $\rho_{n,S},\ \tau_n,\ \tilde{r}_n$ are preferable, otherwise. Of course, the range "close``\ is negotiable.  It was also shown  that amongst the rank coefficients    $\tau_n,\ \rho_{n,S}$ and $\tilde{r}_n$ the coefficient  $\tilde{r}_n$  has the smallest variances when the association between $X$ and $Y$ is not "close``\ to $\pm 1$. Some   disadvantages of $\rho/\rho_n$ were  also highlight in Stepanov (2025b). Thus, the coefficient $\rho_n$ is very sensitive to contamination of outliers, while $\rho_{n,S},\ \tau_n$ and $\tilde{r}_n$ are rather robust. Another disadvantage of $\rho$ and, correspondingly of $\rho_n$, is that $\rho$ cannot be defined  if  the  second moments do not exist. Observe that $\rho_S,\ \tau$ and $r$ exist for any absolutely continuous non-degenerate distribution function.

In the present work, we continue to study the correlation coefficients. In Section~2, we introduce a new rank correlation coefficient $r_n$, which is a sample analogue of $r$, but it is simpler than $\tilde{r}_n$. We study the properties of $r_n$ and compare $r_n$  with $\rho_n,\ \rho_{n,S}$ and $\tau_n$. In particular, we  find out what the asymptotic behaviors of $Var(\tau_n)$ and $Var(r_n)$ look like in the case of an arbitrary absolutely continuous $F$. We also conduct simulation experiments  and discuss the properties of the sample variances of correlation coefficients. All the proofs for the results of this work are presented in Section~3 (Appendix).

\section{New Rank Correlation Coefficient, Results and Examples}

Before we introduce a new rank correlation coefficient,  we compare the rank coefficients $\tau_n$ and $\rho_{n,S}$. Although the coefficient $\rho_{n,S}$ works conscientiously, it is complexly organized.  It is extremely difficult to find its variance analytically for further use in any theoretical analysis.  In addition, there is a redundancy  in the definition of $\rho_{n,S}$. Indeed, $I_{ji}=1-I_{ij}$. The coefficient $\tau _n$ is simply and nicely constructed, so, it can be  studied  theoretically. However, it has one small drawback. In the sum $\sum_{i=2}^n\sum_{j=1}^{i-1}I_{ji}$, the indicators are taken with equal weights (ones). It is known, see, for example Walsh (1969), that  the association between distant order statistics is weaker than between neighboring ones. The same is true for the concomitants. Our basic idea is the following. In the sum $\sum_{i=2}^n\sum_{j=1}^{i-1}I_{ji}$, we substitute $I_{ji}$ with $(n-i+j)I_{ji}$. Then, the indicators that reflect a strong dependence between neighboring concomitants will be multiplied by the largest weights $n-1$, the weights  decrease as the distance between indices in the indicators increases, so that the indicators that reflect the weakest dependence between distant concomitants will be multiplied by the smallest weights $1$. In our opinion, now  $r_n$ shows the true dependence rate between $X$ and $Y$ most accurately. In turn, this accuracy is reflected by the corresponding variance of the correlation coefficient.

\noindent {\bf New Rank Correlation Coefficient}\ Let $T_n=\sum_{i=2}^n\sum_{j=1}^{i-1}(n-i+j)I_{ji}$. By the above argument, we define a new rank correlation coefficient $r_n$ as 
$$
r_n=\frac{12T_n}{n(n-1)(2n-1)}- 1=\frac{12\sum_{i=2}\sum_{j=1}^{i-1}(n-i+j)I_{ji}}{n(n-1)(2n-1)}-1.
$$
Since $\sum_{i=2}^n\sum_{j=1}^{i-1}(n-i+j)=\frac{n(n-1)(2n-1)}{6}$, it is easily seen that $-1\leq r_n\leq 1\ a.s.$ Some other properties of $r_n$ are presented in Theorems~\ref{theorem2.1}, \ref{theorem2.2} in Subsection~2.1.

\subsection{Results}

\begin{theorem}\label{theorem2.1} The expected value of $r_n$ has the form
$$
Er_n=\left(1-\frac{3}{2n-1}\right)\cdot r-\frac{3}{2n-1}+\frac{12E(F(X,Y))}{2n-1}.
$$
\end{theorem}
It follows from Theorem~\ref{theorem2.1} that $Er_n\rightarrow r\ (n\rightarrow \infty)$.  Let us further use the  designations
$$
d_n=dx_1dy_1\ldots dx_ndy_n\ (n\geq 1),\quad f_i=f(x_i,y_i),
$$
$$
H_i=H(x_i),\quad G_i=G(y_i),\quad F_i=F(x_i,y_i),\quad \bar{F}_i=\bar{F}(x_i,y_i),
$$
$$
A(x,y)=\int_{x_1\leq x,\ y_1\leq y}f_1H_1d_1-\int_{x_1\geq  x,\ y_1\geq  y}f_1H_1d_1,
$$
$$
B(x)=\int_{x\leq x_1,\ y_1\in \mathbb{R}}f_1F(x,y_1)d_1.
$$
Let also
$$
Q_1=-4\int_{x_1\leq x_2}f_1f_2F_1F(x_1,y_2)d_2,
$$
\begin{eqnarray*}
Q_2&=&4\int_{x_1,y_1\in \mathbb{R}}f_1H_1[F_1-\bar{F}_1]B(x_1)d_1\\
&+&2\int_{x_1\geq x_2,\ y_1\geq y_2}f_1f_2H_1B(x_2)d_2-2\int_{x_1\leq x_2,\ y_1\leq y_2}f_1f_2H_1B(x_2)d_2,
\end{eqnarray*}
$$
Q_3=\int_{x_1,y_1 \in \mathbb{R}}f_1[F_1+\bar{F}_1]A(x_1,y_1)d_1
$$
and
\begin{eqnarray*}
Q_4=\int_{x_1,y_1 \in \mathbb{R}} f_1A^2(x_1,y_1)d_1+ 2\int_{x_1,y_1 \in \mathbb{R}} f_1H_1[\bar{F}_1-F_1]A(x_1,y_1)d_1.
\end{eqnarray*}
Below we present asymptotic results for the variances of   $r_n$ and $\tau_n$.

\begin{theorem}\label{theorem2.2} For any absolutely continuous distribution function $F$, the asymptotic variance of $r_n$ is as follows
\begin{eqnarray*}
Var(r_n)&=&\frac{36}{n}\left( E[(1+H(X))\bar{F}(X,Y)+(1-H(X))F(X,Y)]^2-\left(\frac{\tau+1}{2}-\frac{\rho_S+1}{4}\right)^2\right)\\
&+& \frac{36}{n}(Q_1+Q_2+Q_3+Q_4)+O\left(\frac{1}{n^2}\right).
\end{eqnarray*}
\end{theorem}

\begin{theorem}\label{theorem2.3} For any absolutely continuous distribution function $F$, the asymptotic variance of $\tau_n$ is as follows
$$
Var(\tau_n)=\frac{16(E[F(X,Y)+\bar{F}(X,Y)]^2-4[EF(X,Y)]^2)}{n}+O\left(\frac{1}{n^2}\right).
$$
\end{theorem}
\begin{remark}\label{remark2.1}
Since $E[\bar{F}(X,Y)+F(X,Y)]^2\geq [E\bar{F}(X,Y)+EF(X,Y)]^2$, we obtain  that
\begin{eqnarray*}
E(\bar{F}(X,Y)+F(X,Y))^2&-&4(EF(X,Y))^2\\
&\geq& [E\bar{F}(X,Y)-EF(X,Y)][E\bar{F}(X,Y)+3EF(X,Y)]=0.
\end{eqnarray*}
That way, the main term of $Var(\tau_n)$ is non-negative. Unfortunately, we cannot show the same for $Var(r_n)$.
\end{remark}
Theorems~\ref{theorem2.2}, \ref{theorem2.3} allow us to compare $Var(\tau_n)$ and $Var(r_n)$ for large $n$ theoretically and decide which correlation coefficient is better  for a  given particular distribution. The later can be shown for all or  for some parametric values.  Further discussion of the properties of correlation coefficients is carried out using examples.

\subsection{Examples}
In this subsection, we discuss three illustrative examples and conduct simulation experiments. 

\noindent {\bf Example~2.1}\ Let $F_t$ be a Farlie--Gumbel--Morgenstern copula with
$$
F_t(x,y)=xy+t(x-x^2)(y-y^2)\quad (-1\leq t\leq 1,\ x,y\in[0,1]).
$$
We have
$$
\rho(t)=\rho_S(t)=\frac{6t}{18},\quad \tau(t)=\frac{4t}{18},\quad r(t)=\frac{3t}{18}
$$
and, then
$$
\left(\frac{\tau(t)+1}{2}-\frac{\rho_S(t)+1}{4}\right)^2=1/16+t/72+t^2/1296.
$$
After tedious computation, we have found that 
$$
E[(1+H(X))\bar{F}(X,Y)+(1-H(X))F(X,Y)]^2=\frac{41}{180}+\frac{t}{12}+\frac{t^2}{180},
$$
$$
E(\bar{F}(X,Y)+F(X,Y))^2=5/18+t/9+t^2/150, 
$$
$$
(EF(X,Y))^2=4\left(\frac{1}{4}+\frac{t}{18}\right)^2,
$$
$$
Q_1=-1/12-2t/45-t^2/180,\quad Q_2=1/30+23t/1080+11t^2/3240,
$$
$$
Q_3=-1/12-2t/45-t^2/180,\quad Q_4=-1/40-t/540+t^2/540.
$$
Observe that for this distribution $Q_1=Q_3$. It follows that 
$$
Var(r_n(t))=\frac{1}{n}\left(\frac{1}{4}-\frac{7}{180}t^2\right)+O\left(\frac{1}{n^2}\right)
$$
and
$$
Var(\tau_n(t))=\frac{1}{n}\left(\frac{4}{9}-\frac{46}{2025}t^2\right)+O\left(\frac{1}{n^2}\right).
$$
Obviously, $Var(r_n(t))<Var(\tau_n(t))$ for   large $n$ and all $t\in[-1,1]$. We   generated one thousand times the correlation coefficients $\rho_n, \rho_{n,S}, \tau_n$ and $r_n$ for samples $X_1,\ldots, X_{1000}$ and $Y_1,\ldots, Y_{1000}$, computed their sample variances and, for visual comparison, presented them and the main terms of the variances of $\tau_n$ and $r_n$ in   Table~2.1 for different positive $t$. For negative $t$ the situation is similar and we do not present the corresponding  results here.

\begin{table}[ht]
\begin{center}
{\em Table~2.1 Sample variances and variances of correlation coefficients for different $t$.}

\vspace*{1ex}
\begin{tabular}{|ccccccc|}
\hline
                   &$t=0.01$              &$t=0.30$     &    $t=0.5$        & $t=0.70$         & $t=0.99$   &     \\
\hline
$S^2_{\rho_{n}}$   &$9.727\cdot 10^{-4}$ &$9.652\cdot 10^{-4}$&$9.488\cdot 10^{-4}$&$8.510\cdot 10^{-4}$&  $7.109\cdot 10^{-4}$&\\
\hline
$S^2_{\rho_{n,S}}$ &$9.708\cdot 10^{-4}$ &$9.631\cdot 10^{-4}$&$9.498\cdot 10^{-4}$&  $8.520\cdot 10^{-4}$&  $7.108\cdot 10^{-4}$&\\
\hline
$S^2_{\tau_n}$     &$4.326\cdot 10^{-4}$ &$4.309\cdot 10^{-4}$&$4.279\cdot 10^{-4}$&  $3.898\cdot 10^{-4}$&  $3.340\cdot 10^{-4}$&\\
\hline
$Var(\tau_n)$      &$4.444\cdot 10^{-4}$ &$4.424\cdot 10^{-4}$&$4.388\cdot 10^{-4}$&$4.333\cdot 10^{-4}$&$4.221\cdot 10^{-4}$&\\
\hline
$S^2_{r_n}$        &$2.432\cdot 10^{-4}$ &$2.435\cdot 10^{-4}$&$2.437\cdot 10^{-4}$&  $2.260\cdot 10^{-4}$&  $1.995\cdot 10^{-4}$& \\
\hline
$Var(r_n)$         &$2.500\cdot 10^{-4}$ &$2.465\cdot 10^{-4}$&$2.403\cdot 10^{-4}$ & $2.309\cdot 10^{-4}$ & $2.119\cdot 10^{-4}$ &\\
\hline
\end{tabular}
\end{center}
\end{table}
The simulation results show a correspondence between $S^2_{\tau_n},\ S^2_{r_n}$, on the one hand, and $Var(\tau_n),\ Var(r_n)$, on the other hand. With unimportant exception for $S^2_{\rho_{n}},\ S^2_{\rho_{n,S}}$ (when $t=0.5,\ 0.7$) we have 
$$
S^2_{r_n}\leq S^2_{\tau_n}\leq  S^2_{\rho_{n,S}}\leq S^2_{\rho_{n}}.
$$
We see that the performance of $r_n$ is much better than the performance of other correlation coefficients. 

\noindent {\bf Example~2.2}\ Let
$$
F_t(x,y)=\frac{1}{2\pi\sqrt{1-t^2}}\int_{-\infty}^x\int_{-\infty}^y e^{-\frac{u^2-2t uv+v^2}{2(1-t^2)}}dudv\quad (x,y\in \mathbb{R},\ -1<t<1)
$$
be a bivariate normal distribution. A deep comparative analysis of $\rho_n, \rho_{n,S}, \tau_n$ was done in Xu et al. (2013). Later a thorough simulation analysis of  $\rho_n, \rho_{n,S}, \tau_n$ and $\tilde{r}_n$ was also carried out in Stepanov (2025b).  In this work, we briefly repeat the last analysis  with using $r_n$ instead of $\tilde{r}_n$.
It should be mentioned that  it was shown in Hotelling (1953) and Esscher (1924) by indirect methods that 
\begin{eqnarray}
Var (\rho_n)&=& \frac{(1-t^2)^2}{n}+O\left(\frac{1}{n^2}\right)\quad (n\rightarrow \infty),\label{2.1}\\
Var (\tau_n)&=&\frac{4}{n}\left(\frac{1}{9}-\frac{4\cdot \arcsin^2(t/2)}{\pi^2}\right)+O\left(\frac{1}{n^2}\right)\quad (n\rightarrow \infty).\label{2.2}
\end{eqnarray}
Observe that finding  $Var(r_n)$ analytically by Theorems~\ref{theorem2.2} in the normal case is not easy. It can be done in our (or someone else) future works. 
It is known, see, for example, Xu et al. (2013) or Stepanov (2025b), that
$$
\rho(t)=t,\quad \rho_S(t)=\frac{6}{\pi}\arcsin(t/2),\quad \tau(t)=\frac{2}{\pi}\arcsin(t)
$$
and
$$
r(t)=\frac{3}{\pi}\cdot (\arcsin(t)-\arcsin(t/2)).
$$
For this work, in turn,  we  generated one thousand times the correlation coefficients $\rho_n, \rho_{n,S}, \tau_n$ and $r_n$ for samples $X_1,\ldots, X_{1000}$ and $Y_1,\ldots, Y_{1000}$, computed their sample variances and presented them  in   Table~2.2 for different positive $t$. For negative $t$ the situation is similar and we do not present the corresponding  results here.
\begin{table}[ht]
\begin{center}
{\em Table~2.2 Sample variances of correlation coefficients for different $t$.}

\vspace*{1ex}
\begin{tabular}{|cccccc|}
\hline
                   & $t=0.05$                & $t=0.30$           & $t=0.70$           & $t=0.99$ &\\
\hline
$S^2_{\rho_{n}}$   & $9.876\cdot 10^{-4}$    & $8.578\cdot 10^{-4}$ & $2.606\cdot 10^{-4}$ & $3.924\cdot 10^{-7}$ &\\
\hline
$S^2_{\rho_{n,S}}$ &  $10.067 \cdot 10^{-4}$ & $8.851\cdot 10^{-4}$ & $3.218\cdot 10^{-4}$ & $7.310\cdot 10^{-7}$ & \\
\hline
$S^2_{\tau_n}$     &  $4.491 \cdot 10^{-4}$  & $ 4.208\cdot 10^{-4}$ & $2.322\cdot 10^{-4}$ & $1.017\cdot 10^{-5}$ &\\
\hline
$S^2_{r_n}$        &  $2.528 \cdot 10^{-4}$  & $2.528\cdot 10^{-4}$ & $1.956\cdot 10^{-4}$ & $1.904\cdot 10^{-5}$ &\\
\hline
\end{tabular}
\end{center}
\end{table}
Table~2.2 tells us  that when the association between $X$ and $Y$ is close to linearity the correlation coefficient  $\rho_n$ works better than any other correlation coefficient. To be precise,  we have $S^2_{\rho_{n}}<S^2_{\rho_{n,S}}<S^2_{\tau_{n}}<S^2_{r_{n}}$ for $t\in W=(-1,-0.730072)\cup (0.730072,1)$. The set $W$ was established in Stepanov (2025b) by the help of the asymptotic relative efficiency/ARE of the  sample variances of correlation coefficients. First, it was obtained theoretically by ARE  that  the coefficient  $\rho_n$ performs asymptotically better than the coefficient  $\tau_n$ when $t\in W$ and, on the contrary, $\tau_n$ performs asymptotically better than $\rho_n$ when $t\in (-1,1)\setminus  W$. By simulation this property was extended for other correlation coefficients.

 \noindent {\bf Example~2.3}\ Let us consider the bivariate Pareto distribution function
$$
F_t(x,y)=1-\frac{1}{(y+1)^t}-\frac{1}{(x+1)^t}+\frac{1}{(x+y+1)^t}\quad (x>0,\ y>0,\ t>0).
$$
We can find that $\rho(t) =\frac{1}{t}\ (t>2)$ and   $\tau(t)=\frac{1}{2t+1}\ (t>0)$. Since the correlation coefficients $\rho_s(t)$ and $r(t)$ can be found  only numerically, we present   in Table~2.3   their values for different  $t$ along with the values of $\rho(t)$ and   $\tau(t)$. 
\begin{table}[ht]
\begin{center}
{\em Table~2.3 Values of correlation coefficients for different $t$.}

\vspace*{1ex}
\begin{tabular}{|ccccccc|}
\hline
  & $t=0.05$ & $t=1$ & $t=2.1$ &  $t=10$ & $t=50$ & $t=100$\\
 \hline
 $\rho$  & -- & -- & 0.4761 &  0.1000 & 0.0200 & 0.0100\\
\hline
$\rho_s$  & 0.6455 & 0.4784 & 0.2839 &  0.0714 & 0.0149 & 0.0075\\
\hline
$\tau$   & 0.9091 & 0.3333 &  0.1923 &  0.0476 & 0.0099 & 0.0050\\
\hline
$r$    & 0.5088 & 0.2608 & 0.1465 &  0.0357 & 0.0074 & 0.0037\\
\hline
\end{tabular}
\end{center}
\end{table}
The simulation experiments show  that $\rho_n(t)$ measures the dependence rate when $t\in(0,2]$ and when its theoretical analogue $\rho(t)$ does not exist. In connection with this, an extension of  $\rho(t)\ (t\in(0,2])$ was proposed in Stepanov (2025b). We again    generated one thousand times the correlation coefficients $\rho_n, \rho_{n,S}, \tau_n$ and $r_n$ for samples $X_1,\ldots, X_{1000}$ and $Y_1,\ldots, Y_{1000}$, computed their sample variances and presented them  in   Table~2.4 for different positive $t$.
\begin{table}[ht!]
\begin{center}
{\em Table~2.4 Sample variances  of correlation coefficients for different $t$.}

\vspace*{1ex}
\begin{tabular}{|cccccccc|}
\hline
                   &  $t=0.05$  & $t=1$ & $t=2.1$ &  $t=10$ & $t=50$ & $t=100$&\\
 \hline
 $S^2_{\rho_{n}}$  & $1.360\cdot 10^{-2}$ & $5.227\cdot 10^{-2}$ & $2.178\cdot 10^{-2}$ & $1.643\cdot 10^{-3}$ & $1.076\cdot 10^{-3}$& $1.136\cdot 10^{-3}$&\\
\hline
$S^2_{\rho_{n,S}}$ & $1.641\cdot 10^{-6}$ & $7.062\cdot 10^{-4}$ & $9.644\cdot 10^{-4}$ & $1.009\cdot 10^{-3}$& $9.914\cdot 10^{-4}$& $1.127\cdot 10^{-3}$&\\
\hline
$S^2_{\tau_n}$     & $1.641\cdot 10^{-5}$ & $3.913\cdot 10^{-4}$ & $4.654\cdot 10^{-4}$ & $4.532\cdot 10^{-4}$ & $4.429\cdot 10^{-4}$& $5.020\cdot 10^{-4}$&\\
\hline
 $S^2_{r_n}$       & $2.976\cdot 10^{-5}$ & $2.702\cdot 10^{-4}$ & $2.837\cdot 10^{-4}$ & $2.569\cdot 10^{-4}$& $2.496\cdot 10^{-4}$& $2.821\cdot 10^{-4}$&\\
\hline
\end{tabular}
\end{center}
\end{table}
We see that  $r_n$ performs better than  other correlation coefficients with exception when $t=0.05$. By further simulation, we discovered that $\rho_{n,S}$ and $\tau_n$ perform a little better than $r_n$ in the interval $[0.0012, 0.07]$.

{\bf Conclusion}\ In the present paper, we have proposed  a new rank correlation coefficient $r_n$, which is a sample analogue of the theoretical correlation coefficient $r$, which, in turn, was previously proposed in Stepanov (2025b). We have discussed the properties of $r_n$ and compared $r_n$ with the known rank Spearman and Kendall correlation coefficients  as well as with the sample Pearson correlation coefficient.  This comparison and simulation experiments   show   that in the case when the association between $X$ and $Y$ is not close to linearity,  $r_n$ performs better than other correlation coefficients. We also have established the asymptotic behavior of $Var(\tau_n)$ and $Var(r_n)$ in the case of arbitrary absolutely continuous distribution.

\section{Appendix}

\begin{gproof}{ of Theorem~\ref{theorem2.1}} 
The joint density function of the vector $[X_{(j)},\ Y_{[j]},\ X_{(i)},\ Y_{[i]}]\ (j<i)$ for   $x_1<x_2,\ y_1,y_2\in \mathbb{R}$ has, see, for example, Balakrishnan and Stepanov (2009), the  form
\begin{equation}\label{3.1}
 f(x_1,y_1,x_2,y_2)=\frac{n!}{(j-1)!(i-j-1)!(n-i)!} f_1f_2H^{j-1}_1(H_2-H_1)^{i-j-1}(1-H_2)^{n-i}. 
\end{equation}
It follows that
\begin{eqnarray*}
ET_n&=&\sum_{i=2}^n\sum_{j=1}^{i-1}(n-i+j)P(Y_{[j]}\leq Y_{[i]})\\
&=&\sum_{i=2}^n\sum_{j=1}^{i-1} \frac{(n-1)n!}{(j-1)!(i-j-1)!(n-i)!} \int_{x_1\leq x_2,\ y_1\leq y_2} f_1f_2H^{j-1}_1(H_2-H_1)^{i-j-1}(1-H_2)^{n-i}d_2\\
&-&\sum_{i=2}^n\sum_{j=1}^{i-1} \frac{(i-j-1)n!}{(j-1)!(i-j-1)!(n-i)!} \int_{x_1\leq x_2,\ y_1\leq y_2} f_1f_2H^{j-1}_1(H_2-H_1)^{i-j-1}(1-H_2)^{n-i}d_2\\
&=&n(n-1)(n-1)\int_{x_1\leq x_2,\ y_1\leq y_2} f_1f_2d_2-n(n-1)(n-2)\int_{x_1\leq x_2,\ y_1\leq y_2} f_1f_2(H_2-H_1)d_2\\
&=&n(n-1)^2EF(X,Y)+ n(n-1)(n-2)[1/6-E(H(X)G(Y))].
\end{eqnarray*}
By making use of the identity $r=\frac{3}{2}\tau-\frac{1}{2}\rho_{S}$, we obtain
$$
Er_n=\left(1-\frac{3}{2n-1}\right)\cdot r-\frac{3}{2n-1}+\frac{12E(F(X,Y))}{2n-1}.
$$
\end{gproof}

It follows from the above proof that 
\begin{eqnarray}\label{3.2}
\nonumber &&(ET_n)^2=O(n^4)+n^6[E(F(X,Y))+1/6-E(H(X)G(Y))]^2-\\
 && \hspace{-1ex}n^5[4(EF(X,Y))^2+10E(F(X,Y))[1/6-E(H(X)G(Y))]+6(1/6-E(H(X)G(Y)))^2].
\end{eqnarray}

\begin{gproof}{ of Theorem~\ref{theorem2.2}} 
Observe that 
$$
Var(r_n)=\frac{144Var(T_n)}{n^2(n-1)^2(2n-1)^2}.
$$
Let us find $ET_n^2$. We have
\begin{eqnarray*}
ET^2_n&=&\sum_{i=2}^n\sum_{j=1}^{i-1}\sum_{k=2}^n\sum_{m=1}^{k-1}(n-k+m)(n-i+j)P(Y_{[j]}\leq  Y_{[i]},Y_{[m]}\leq Y_{[k]})\\
&=&\Sigma_1+\Sigma_2.
\end{eqnarray*}
Here
\begin{eqnarray*}
\Sigma_1=\sum_{j<i<m<k}+\sum_{j<m<i<k}+\sum_{m<j<i<k}+\sum_{m<j<k<i}+\sum_{m<k<j<i}+\sum_{j<m<k<i}
\end{eqnarray*}
and
\begin{eqnarray*}
\Sigma_2&=&\sum_{m<j<k=i}+\sum_{j<m<k=i}+\sum_{j<m=i<k}\\
&+&\sum_{m<k=j<i}+\sum_{m=j<i<k}+\sum_{m=j<k<i}+\sum_{m=j<k=i},
\end{eqnarray*}
where, in turn, for example,  
$$
\sum_{j<i<m<k}=\sum_{k=4}^n\sum_{m=3}^{k-1}\sum_{i=2}^{m-1}\sum_{j=1}^{i-1}(n-k+m)(n-i+j)P(Y_{[j]}\leq  Y_{[i]},Y_{[m]}\leq Y_{[k]})
$$
and 
$$
\sum_{m<j<k=i}=\sum_{i=3}^{n}\sum_{j=2}^{i-1}\sum_{m=1}^{j-1}(n-i+m)(n-i+j)P(Y_{[j]}\leq  Y_{[i]},Y_{[m]}\leq Y_{[i]}).
$$
In the proof of Theorem~\ref{theorem2.2} we present two technical propositions.
\begin{proposition}\label{proposition3.1}
The identity
\begin{eqnarray*}
\Sigma_1&=&(n^6-8n^5)(EF(X,Y))^2+(n^6-15n^5)[1/6-E(H(X)G(Y))]^2\\
&+&2(n^6-11n^5)E(F(X,Y))[1/6-E(H(X)G(Y))]+n^5(Q_1+Q_2)+O(n^4).
\end{eqnarray*}
holds true.
\end{proposition}

\begin{gproof}{ of Proposition~\ref{proposition3.1}} 
The joint density function of the vector 
$$
[X_{(m)},\ Y_{[m]},\ X_{(j)},\ Y_{[j]},\ X_{(i)},\ Y_{[i]},\ X_{(k)},\ Y_{[k]}]\quad (m<j<i<k)
$$ 
for  $x_1\leq x_2\leq x_3\leq x_4,\ y_1,y_2,y_3,y_4\in \mathbb{R}$ has  (see, for example, Balakrishnan and Stepanov (2009)) the  form
\begin{eqnarray}\label{3.3}
\nonumber f(x_1, y_1, x_2, y_2, x_3, y_3, x_4, y_4)&=&\frac{n!}{(m-1)!(m-j-1)!(i-j-1)!(k-i-1)!(n-i)!}\times f_1f_2f_3f_4\\
&\times& H^{m-1}_1H^{j-m-1}_2(H_3-H_2)^{i-j-1}(H_4-H_3)^{k-i-1}(1-H_4)^{n-k}.
\end{eqnarray}
Let us consider $\sum_{j<i<m<k}$. We can write
\begin{eqnarray*}
(n-k+m)(n-i+j)&=&(n-1)^2\\
&-&(n-1)(k-m-1)-(n-1)(i-j-1)+(k-m-1)(i-j-1).
\end{eqnarray*}
Then,
$$
\sum_{j<i<m<k} =S_{(n-1)^2}-S_{(n-1)(k-m-1)}-S_{(n-1)(i-j-1)}+S_{(k-m-1)(i-j-1)},
$$
where, for example,
$$
S_{(n-1)^2}=\sum_{k=4}^n\sum_{m=3}^{k-1}\sum_{i=2}^{m-1}\sum_{j=1}^{i-1}(n-1)^2P(Y_{[j]}\leq  Y_{[i]},Y_{[m]}\leq Y_{[k]})
$$
and
$$
S_{(n-1)(k-m-1)}=\sum_{k=4}^n\sum_{m=3}^{k-1}\sum_{i=2}^{m-1}\sum_{j=1}^{i-1}(n-1)(k-m-1)P(Y_{[j]}\leq  Y_{[i]},Y_{[m]}\leq Y_{[k]}).
$$
We have
\begin{eqnarray*}
S_{(n-1)^2}&=&\sum_{k=4}^n\sum_{m=3}^{k-1}\sum_{i=2}^{m-1}\sum_{j=1}^{i-1}\frac{n!(n-1)^2}{(j-1)!(i-j-1)!(m-i-1)!(k-i-1)!(n-k)!}\\
&& \hspace{-12ex}\times \int_{x_1\leq x_2\leq x_3\leq x_4,\ y_1\leq y_2,\ y_3\leq y_4}f_1f_2f_3f_4H^{j-1}_1H^{i-j-1}_2(H_3-H_2)^{m-i-1}(H_4-H_3)^{k-m-1}(1-H_4)^{n-k}d_4.\\
&=& n(n-1)(n-2)(n-3)(n-1)^2I_{1234,12,34}=(n^6-8n^5)I_{1234;12,34}+O(n^4),
\end{eqnarray*}
where 
$$
I_{1234;12,34}=\int_{x_1\leq x_2\leq x_3\leq x_4,\ y_1\leq y_2,\ y_3\leq y_4}f_1f_2f_3f_4d_4.
$$
Computing the other terms of $\sum_{j<i<m<k}$ and summing the results, we obtain
\begin{eqnarray*}
\sum_{j<i<m<k}&=&(n^6-8n^5)I_{1234;12,34}\\
&+&(n^6-11n^5)[I^{(1)}_{1234;12,34}-I^{(2)}_{1234;12,34}+I^{(3)}_{1234;12,34}-I^{(4)}_{1234;12,34}]+(n^6-15n^5)I_{1234;12,34}^{(4-2)(3-1)},
\end{eqnarray*}
where 
$$
I^{(l)}_{1234;12,34}=\int_{x_1\leq x_2\leq x_3\leq x_4,\ y_1\leq y_2,\ y_3\leq y_4}f_1f_2f_3f_4H_ld_4\quad (l=1,2,3,4)
$$
and
$$
I_{1234;12,34}^{(4-2)(3-1)}=\int_{x_1\leq x_2\leq x_3\leq x_4,\ y_1\leq y_2,\ y_3\leq y_4}f_1f_2f_3f_4(H_4-H_2)(H_3-H_1)d_4.
$$
Similarly, we get
\begin{eqnarray*}
\sum_{j<m<i<k}&=&(n^6-10n^5)I_{1324;12,34}+(n^6-12n^5)(I^{(1)}_{1324;12,34}-I^{(4)}_{1324;12,34})\\
&+& (n^6-13n^5)(I^{(3)}_{1324;12,34}-I^{(2)}_{1324;12,34})+(n^6-15n^5)I_{1324;12,34}^{(4-2)(3-1)},
\end{eqnarray*}
\begin{eqnarray*}
\sum_{m<j<i<k}&=&(n^6-10n^5)I_{3124;12,34}+(n^6-14n^5)(I^{(1)}_{3124;12,34}-I^{(2)}_{3124;12,34})\\
&+&(n^6-11n^5)(I^{(3)}_{3124;12,34}-I^{(4)}_{3124;12,34})+(n^6-15n^5)I_{1324;12,34}^{(4-2)(3-1)},
\end{eqnarray*}
\begin{eqnarray*}
\sum_{m<j<k<i}&=&(n^6-10n^5)I_{3142;12,34}+(n^6-12n^5)(I^{(3)}_{3142;12,34}-I^{(2)}_{3142;12,34})\\
&+&(n^6-13n^5)(I^{(1)}_{3142;12,34}-I^{(4)}_{3142;12,34})+(n^6-15n^5)I_{1324;12,34}^{(4-2)(3-1)},
\end{eqnarray*}
\begin{eqnarray*}
\sum_{m<k<j<i}&=&(n^6-8n^5)I_{3412;12,34}\\
&+&(n^6-11n^5)[I^{(1)}_{3412;12,34}-I^{(2)}_{3412;12,34}I^{(3)}_{3412;12,34}-I^{(4)}_{3412;12,34}]+(n^6-15n^5)I_{1324;12,34}^{(4-2)(3-1)},
\end{eqnarray*}
\begin{eqnarray*}
\sum_{j<m<k<i}&=&(n^6-10n^5)I_{1342;12,34}+(n^6-11n^5)(I^{(1)}_{1342;12,34}-I^{(2)}_{1324;12,34})\\
&+&(n^6-14n^5)(I^{(3)}_{1342;12,34}-I^{(4)}_{1342;12,34})+(n^6-15n^5)I_{1324;12,34}^{(4-2)(3-1)}.
\end{eqnarray*}
That way,
$$
\Sigma_1= (n^6-8n^5)K_1+n^5Q_1+(n^6-11n^5)K_2+n^5Q_2+(n^6-15n^5)J+O(n^4),\\
$$
where
$$
K_1=I_{1234;12,34}+I_{1324;12,34}+I_{3124;12,34}+I_{3142;12,34}+I_{3412;12,34}+I_{1342;12,34},
$$
$$
Q_1=-2(I_{1324;12,34}+I_{3124;12,34}+I_{3142;12,34}+I_{1342;12,34})
$$
and
$$
K_2=\sum (I^{(1)}_{l_1l_2l_3l_4;12,34}-I^{(2)}_{l_1l_2l_3l_4;12,34}+I^{(3)}_{l_1l_2l_3l_4;12,34}-I^{(4)}_{l_1l_2l_3l_4;12,34}).
$$
There are 6 terms in $K_2$, where  the set $[l_1,l_2,l_3,l_4]$ is equal  to one of the following six sets $1234,\ 1324,\ 3124,\ 3142,\ 3412,\ 1342$. The other terms of $\Sigma_1$ have the forms:
\begin{eqnarray*}
Q_2&=& -I_{1324;12,34}^{(1)}+2I_{1324;12,34}^{(2)}-2I_{1324;12,34}^{(3)}+I_{1324;12,34}^{(4)}-3I_{3124;12,34}^{(1)}+3I_{3124;12,34}^{(2)}\\
&-& 2I_{3142;12,34}^{(1)}+I_{3142;12,34}^{(2)}-I_{3142;12,34}^{(3)}+2I_{3142;12,34}^{(4)}-3I_{1342;12,34}^{(3)}+3I_{1342;12,34}^{(4)}
\end{eqnarray*}
and
$$
J=I_{1234;12,34}^{(4-2)(3-1)}+I_{1342;12,34}^{(4-2)(3-1)}+I_{3124;12,34}^{(4-2)(3-1)}+I_{3142;12,34}^{(4-2)(3-1)}+I_{3412;12,34}^{(4-2)(3-1)}+I_{1342;12,34}^{(4-2)(3-1)}.
$$
Observe that
$$
K_1=\int_{x_1\leq x_2,\  x_3\leq x_4,\ y_1\leq y_2,\ y_3\leq y_4}f_1f_2f_3f_4d_4=(EF(X,Y))^2.
$$
Since
$$
I_{1324;12,34}+I_{3124;12,34}=I_{3142;12,34}+I_{1342;12,34}=\int_{x_1\leq x_2}f_1f_2F_1F(x_1,y_2)d_2,
$$
we get
$$
Q_1=-4\int_{x_1\leq x_2}f_1f_2F_1F(x_1,y_2)d_2.
$$
We also have
$$
K_2= 2E(F(X,Y))[1/6-E(H(X)G(Y))],
$$
Let us consider $Q_2$. Let us take from $Q_2$ the terms $2I_{1324;12,34}^{(2)}+3I_{3124;12,34}^{(2)}+I_{3142;12,34}^{(2)}$. We have
\begin{eqnarray*}
2I_{1324;12,34}^{(2)}+3I_{3124;12,34}^{(2)}+I_{3142;12,34}^{(2)}&=&2(I_{1324;12,34}^{(2)}+I_{3124;12,34}^{(2)})+(I_{3124;12,34}^{(2)}+I_{3142;12,34}^{(2)})\\
&& \hspace{-6ex}=2\int_{x_1,y_1\in \mathbb{R}}f_1H_1F_1B(x_1)d_1+\int_{x_1\geq x_1,\ y_1\geq y_2}f_1f_2H_1B(x_2)d_2.
\end{eqnarray*}
Considering the other terms of $Q_2$ 
$$
-I_{1324;12,34}^{(1)}-3I_{3124;12,34}^{(1)}- 2I_{3142;12,34}^{(1)},
$$
$$
-2I_{1324;12,34}^{(3)}-I_{3142;12,34}^{(3)}-3I_{1342;12,34}^{(3)},
$$
$$
I_{1324;12,34}^{(4)}2I_{3142;12,34}^{(4)}3I_{1342;12,34}^{(4)}
$$
and summing up the results, we find that
\begin{eqnarray*}
Q_2&=&4\int_{x_1,y_1\in \mathbb{R}}f_1H_1[F_1-\bar{F}_1]B(x_1)d_1\\
&+&2\int_{x_1\geq y_2,\ y_1\geq y_2}f_1f_2H_1B(x_2)d_2-2\int_{x_1\leq x_2,\ y_1\leq y_2}f_1f_2H_1B(x_2)d_2.
\end{eqnarray*}
It is easily seen that
$$
J=[E(H(X)F(X,Y))-E(H(X)\bar{F}(X,Y))]^2=[1/6-E(H(X)G(Y))]^2.
$$
The result of Proposition~\ref{proposition3.1} readily follows.
\end{gproof}

\begin{proposition}\label{proposition3.2}
The equality
\begin{eqnarray*}
\Sigma_2\\
&&  \hspace{-8ex}  =n^5[E(F(X,Y)+\bar{F}(X,Y))^2+2E(H(X)(\bar{F}^2(X,Y)-F^2(X,Y)))+E[H(X)(\bar{F}(X,Y)-F(X,Y))]^2\\
&+&Q_3+Q_4]+O(n^4)
\end{eqnarray*}
holds true.
\end{proposition}

\begin{gproof}{ of Proposition~\ref{proposition3.2}} The joint density function of the vector 
$$
X_{(j)},\ Y_{[j]},\ X_{(i)},\ Y_{[i]},\ X_{(k)},\ Y_{[k]}\quad  (j<i<k)
$$ 
for  $x_1\leq x_2\leq x_3,\ y_1,y_2,y_3\in \mathbb{R}$ has  (see, for example, Balakrishnan and Stepanov (2009)) the  form
\begin{eqnarray}\label{3.4}
\nonumber f(x_1, y_1, x_2, y_2, x_3, y_3)&=&\frac{n!}{(j-1)!(i-j-1)!(k-i-1)!(n-i)!}\\
&\times& f_1f_2f_3H^{j-1}_1(H_2-H_1)^{i-j-1}(H_3-H_2)^{k-i-1}(1-H_3)^{n-k}.
\end{eqnarray}
Let us consider $\sum_{m<j<k=i}$. We have
\begin{eqnarray*}
\sum_{m<j<k=i}&=& S_{(n-1)(n-2)}-S_{(n-1)(i-j-1)}-S_{(n-1)(j-m-1)}\\
&-&S_{(n-3)(i-j-1)}+S_{(i-j-1)(i-j-2)}+S_{(i-j-1)(j-m-1)},
\end{eqnarray*}
where, for example, 
$$
S_{(n-1)(i-j-1)}=\sum_{i=3}^n\sum_{j=2}^{i-1}\sum_{m=1}^{j-1}(n-1)(i-j-1)P(Y_{[j]}\leq  Y_{[i]},Y_{[m]}\leq Y_{[i]}).
$$
By (\ref{3.4}), we have
\begin{eqnarray*}
S_{(n-1)(i-j-1)}&=&\sum_{i=4}^{n}\sum_{j=2}^{i-2}\sum_{m=1}^{j-1}\frac{n!(n-1)(i-j-1)}{(m-1)!(j-m-1)!(i-j-1)(n-i)!}\\
&& \hspace{-6ex} \times\int_{x_1\leq x_2\leq x_3,\ y_1\leq y_3,\ y_2\leq y_3} f_1f_2f_3H^{m-1}_1(H_2-H_1)^{j-m-1}(H_3-H_2)^{i-m-1}(1-H_3)^{n-i}d_3\\
&=& n^5\int_{x_1\leq x_2\leq x_3,\ y_1\leq y_3,\ y_2\leq y_3} f_1f_2f_3 (H_3-H_2)d_3+O(n^4).
\end{eqnarray*}
In the same way,
\begin{eqnarray*}
S_{(n-1)(n-2)}&=&\sum_{i=3}^{n}\sum_{j=2}^{i-1}\sum_{m=1}^{j-1}\frac{n!(n-1)(n-2)}{(m-1)!(j-m-1)!(i-j-1)(n-i)!}\\
&&\hspace{-6ex} \times\int_{x_1\leq x_2\leq x_3,\ y_1\leq y_3,\ y_2\leq y_3} f_1f_2f_3H^{m-1}_1(H_2-H_1)^{j-m-1}(H_3-H_2)^{i-m-1}(1-H_3)^{n-i}d_3\\
&=& n^5\int_{x_1\leq x_2\leq x_3,\ y_1\leq y_3,\ y_2\leq y_3} f_1f_2f_3 d_3+O(n^4).
\end{eqnarray*}
Computing the other terms of $\sum_{m<j<k=i}$ and summing the results, we obtain
$$
\sum_{m<j<k=i}=n^5\left(I_{123;13,23}+I_{123;13,23}^{(1)}+I_{123;13,23}^{(2)}-2I_{123;13,23}^{(3)}+I^{(3-2)(3-1)}_{123;13,23}\right)+O(n^4),
$$
where
$$
I_{123;13,23}=\int_{x_1\leq x_2\leq x_3,\ y_1\leq y_3,\ y_2\leq y_3} f_1f_2f_3 d_3,
$$
$$
I_{123;13,23}^{(l)}=\int_{x_1\leq x_2\leq x_3,\ y_1\leq y_3,\ y_2\leq y_3} f_1f_2f_3 H_ld_3\quad (l=1,2,3),
$$
and
$$
I^{(3-2)(3-1)}_{123;13,23}=\int_{x_1\leq x_2\leq x_3,\ y_1\leq y_3,\ y_2\leq y_3} f_1f_2f_3 (H_3-H_2)(H_3-H_1)d_3.
$$
Similarly,
$$
\sum_{j<m<k=i}=n^5\left(I_{213;13,23}+I_{213;13,23}^{(1)}+I_{213;13,23}^{(2)}-2I_{213;13,23}^{(3)}+I^{(3-2)(3-1)}_{213;13,23}\right)+O(n^4),
$$
$$
\sum_{j<m=i<k}=n^5\left(I_{123;123}+I_{123;123}^{(1)}-I_{123;123}^{(3)}+I^{(3-2)(2-1)}_{123;123}\right)+O(n^4),
$$
where, for example, 
$$
I_{123;123}=\int_{x_1\leq x_2\leq x_3,\ y_1\leq y_2\leq y_3} f_1f_2f_3 d_3.
$$
We also obtain
$$
\sum_{m<k=j<i}=n^5\left(I_{123;123}+I_{123;123}^{(1)}-I_{123;123}^{(3)}+I^{(3-2)(2-1)}_{123;123}\right)+O(n^4),
$$
$$
\sum_{m=j<i<k}=n^5\left(I_{213;12,13}+2I_{213;12,13}^{(1)}-I_{213;12,13}^{(2)}-I_{213;12,13}^{(3)}+I^{(3-1)(2-1)}_{213;12,13}\right)+O(n^4),
$$
\begin{eqnarray*}
\sum_{m=j<i<k}=n^5\left(I_{132;12,13}+2I_{132;12,13}^{(1)}-I_{132;12,13}^{(2)}-I_{132;12,13}^{(3)}+I^{(3-1)(2-1)}_{132;12,13}\right)+O(n^4),
\end{eqnarray*}
By (\ref{3.1}), we can show that $\sum_{m=j<k=i}=O(n^4)$. We have
\begin{eqnarray*}
L_1&=& \left(I_{123;13,23}+I_{213;13,23}\right)+2I_{123;123}+\left(I_{123;12,13}+I_{132;12,13}\right)\\
&=& E(F^2(X,Y))+2E(F(X,Y)\bar{F}(X,Y))+E(\bar{F}^2(X,Y))=E(F(X,Y)+\bar{F}(X,Y))^2.
\end{eqnarray*}
It should be noted that
\begin{eqnarray*}
I_{123,13,23}^{(1)}+I_{123,13,23}^{(2)}&-&2I_{123,13,23}^{(3)}+I_{213,13,23}^{(1)}+I_{213,13,23}^{(2)}-2I_{213,13,23}^{(3)}\\
&=& 2\int_{x_1\leq x_2,\ y_1\leq y_2}f_1f_2H_1F_2d_2-E(H(X)F^2(X,Y)),
\end{eqnarray*}
$$
2\left(I_{123,123}^{(1)}-I_{123,123}^{(3)}\right)=2\left(\int_{x_1\leq x_2\ y_1\leq y_2} f_1f_2H_1\bar{F}_2-\int_{x_1\geq  x_2\ y_1\geq  y_2} f_1f_2H_1F_2\right),
$$
\begin{eqnarray*}
2I_{123,12,13}^{(1)}-I_{123,12,13}^{(2)}&-&I_{123,12,13}^{(3)}+2I_{132,12,13}^{(1)}-I_{132,12,13}^{(2)}-I_{132,12,13}^{(3)}\\
&=& 2\left(E(H(X)\bar{F}^2(X,Y))-\int_{x_1\geq  x_2,\ y_1\geq  y_2}f_1f_2H_1\bar{F}_2d_2\right).
\end{eqnarray*}
Summing up the last three equalities and denoting the sum by $L_2$, we get
$$
L_2=2(E(H(X)(\bar{F}^2(X,Y)-F^2(X,Y)))+Q_3).
$$
Observe that
\begin{eqnarray*}
I_{123,13,23}^{(3-2)(3-1)}+I_{123,13,23}^{(3-2)(3-1)}&=&-2\int_{x_1\leq x_2,\ y_1\leq y_2}f_1f_2H_1H_2F_2d_2\\
+2\int_{x_1\geq x_2,\ y_1\leq y_2}f_1f_2H_1H_2\bar{F}(x_1,y_2)d_2&+&2\int_{x_1\leq x_2,\ y_1\leq y_2}f_1f_2H_1H_2\bar{F}(_2d_2+E(H^2(X)F^2(X,Y)),
\end{eqnarray*}
\begin{eqnarray*}
I_{123,12,13}^{(3-1)(2-1)}+I_{132,12,13}^{(3-1)(2-1)}&=&-2\int_{x_1\leq x_2,\ y_1\leq y_2}f_1f_2H_1H_2\bar{F}_2\\
+2\int_{x_1\leq x_2,\ y_1\geq  y_2}f_1f_2H_1H_2F(x_1,y_2)&+&2\int_{x_1\geq x_2,\ y_1\geq  y_2}f_1f_2H_1H_2F_2+E(H^2(X)\bar{F}^2(X,Y)),
\end{eqnarray*}
\begin{eqnarray*}
2I_{123,123}^{(3-1)(2-1)}&=&-2E(H^2(X)F(X,Y)\bar{F}(X,Y))\\
&-&2\int_{x_1\leq x_2,\ y_1\leq y_2}f_1f_2H_1H_2F_2d_2+2\int_{x_1\leq x_2,\ y_1\leq y_2}f_1f_2H_1H_2\bar{F}_2d_2\\
&+& 2\int_{x_1\geq  x_2,\ y_1\geq  y_2}f_1f_2H_1H_2F(x_1,y_2)d_2+ 2\int_{x_1\leq   x_2,\ y_1\leq   y_2}f_1f_2H_1H_2F(x_1,y_2)d_2.
\end{eqnarray*}
Summing up the last three equalities and denoting the sum by $L_3$, we get
$$
L_3=E(H^2(X)(\bar{F}^2(X,Y)-F^2(X,Y)))+Q_4.
$$
We have
\begin{eqnarray*}
\Sigma_2&=&L_1+L_2+L_3\\
&=&n^5[E(F(X,Y)+\bar{F}(X,Y))^2+2E(H(X)(\bar{F}^2(X,Y)-F^2(X,Y)))\\
&+&E(H^2(X)(\bar{F}^2(X,Y)-F^2(X,Y)))+Q_3+Q_4]+O(n^4)
\end{eqnarray*}
\end{gproof}
It follows from Propositions~\ref{proposition3.1}, \ref{proposition3.2} that
\begin{eqnarray*}
ET^2_n&=&n^6[EF(X,Y)+1/6-E(H(X)G(Y))]^2\\
 &-& n^5[8(EF(X,Y))^2+22E(F(X,Y))[1/6-E(H(X)G(Y))]+15(1/6-E(H(X)G(Y)))^2]\\
 &&  \hspace{-7ex}+n^5[E(F(X,Y)+\bar{F}(X,Y))^2+2E(H(X)(\bar{F}^2(X,Y)-F^2(X,Y)))+E[H(X)(\bar{F}(X,Y)-F(X,Y))]^2\\
 &+& n^5(Q_1+Q_2+Q_3+Q_4)+O(n^4).
\end{eqnarray*}
The   last identity and (\ref{3.2}) imply that
\begin{eqnarray*}
&& \hspace{-6ex} Var(T_n)\\
 &&\hspace{-6ex} =n^5[E(F(X,Y)+\bar{F}(X,Y))^2+2E(H(X)(\bar{F}^2(X,Y)-F^2(X,Y)))+E[H(X)(\bar{F}(X,Y)-F(X,Y))]^2\\ 
&-& n^5[4(EF(X,Y))^2+12E(F(X,Y))[1/6-E(H(X)G(Y))]+9(1/6-E(H(X)G(Y)))^2]\\
&+& n^5(Q_1+Q_2+Q_3+Q_4)+O(n^4).
\end{eqnarray*}
Finally, observe that
\begin{eqnarray*}
 &&\hspace{-6ex} E(F(X,Y)+\bar{F}(X,Y))^2+2E(H(X)(\bar{F}^2(X,Y)-F^2(X,Y)))+E[H(X)(\bar{F}(X,Y)-F(X,Y))]^2\\ 
&=& E[(1+H(X))\bar{F}(X,Y)+(1-H(X))F(X,Y)]^2
\end{eqnarray*}
and
\begin{eqnarray*}
4(EF(X,Y))^2+12E(F(X,Y))[1/6-E(H(X)G(Y))]&+&9(1/6-E(H(X)G(Y)))^2\\
&=&\left(\frac{\tau+1}{2}-\frac{\rho_S+1}{4}\right)^2.
\end{eqnarray*}

\end{gproof}

\begin{gproof}{ of Theorem~\ref{theorem2.3}} 
The variance of $\tau_n$ was estimated in Stepanov (2025a). Obviously,
$$
Var (\tau_n)=\frac{16}{n^2(n-1)^2} Var\left(\sum_{i=2}^n \sum_{j=1}^{i-1} I_{ji}\right).
$$
Some technical results for $\sum_{i=2}^n \sum_{j=1}^{i-1} I_{ji}$ were  obtained in that paper:
$$
E\left(\sum_{i=2}^n \sum_{j=1}^{i-1} I_{ji}\right)=\int_{x_1\leq x_2,\ y_1\leq y_2}f_1f_2d_2=n(n-1)E(F(X,Y))
$$
and
$$
E\left(\sum_{i=2}^n \sum_{j=1}^{i-1} I_{ji} \right)^2=Q_5+Q_6+Q_7+Q_8,
$$
where
\begin{eqnarray*}
Q_5&=&n(n-1)(n-2)(n-3)\int_{x_1\leq x_2, x_3\leq x_4, y_1\leq y_2, y_3\leq y_4}f_1f_2f_3f_4d_4\\
&=&n(n-1)(n-2)(n-3)E(F(X,Y))^2=(n^4-6n^3)(EF(X,Y))^2+O(n^2),
\end{eqnarray*}
\begin{eqnarray*}
Q_6&=&n(n-1)(n-2)\int_{x_1\leq x_2,\ x_1\leq x_3,\ y_1\leq y_2,\ y_1\leq y_3}f_1f_2f_3d_3\\
&=&n(n-1)(n-2)E(\bar{F}(X,Y))^2=n^3E(\bar{F}(X,Y))^2+O(n^2),
\end{eqnarray*}
\begin{eqnarray*}
Q_7&=&2n(n-1)(n-2)\int_{x_1\leq x_2\leq x_3,\ y_1\leq y_2\leq y_3}f_1f_2f_3d_3\\
&=&2n(n-1)(n-2)E(\bar{F}(X,Y)F(X,Y))=2n^3E(\bar{F}(X,Y)F(X,Y))+O(n^2)
\end{eqnarray*}
and
\begin{eqnarray*}
Q_8&=&n(n-1)(n-2)\int_{x_1\leq x_3,\ x_2\leq x_3,\ y_1\leq y_3,\ y_2\leq y_3}f_1f_2f_3d_3\\
&=&n(n-1)(n-2)E(F(X,Y))^2=n^3E(F(X,Y))^2+O(n^2),
\end{eqnarray*}
It follows that
$$
Var(\tau_n)=\frac{16[E(\bar{F}(X,Y)+F(X,Y))^2-4(EF(X,Y))^2]}{n}+O\left(\frac{1}{n^2}\right).
$$
\end{gproof}

\section*{References}
\begin{description} 
{\small

\item Bairamov, I., Stepanov, A. (2010).\ Numbers of near-maxima for the bivariate case,  {\it Statistics  $\&$ Probability Letters}, {\bf 80}, 196--205.

\item Balakrishnan, N., Lai, C.D. (2009).\ {\it Continuous Bivariate Distributions}, Second edition, Springer.
 
\item Balakrishnan, N., Stepanov, A. (2015).\ Limit results for concomitants of order statistics, {\it Metrika}, {\bf 78}, 385--397.

\item  Bhattacharya, B.B. (1974).\ Convergence of sample paths of normalized sums of induced order statistics, {\it Ann. Statist. }, {\bf 2}, 1034--1039.


\item Daniels, H. E. (1950).\ Rank correlation and population models, {\it Journal of the Royal Statistical Society}, Ser. B, {\bf 12} (2),
171--191.

\item David, H.A., Galambos, J. (1974).\ The asymptotic theory of concomitants of order statistics, {\it J. Appl. Probab.}, {\bf 11}, 762--770.

\item David, H.A., Nagaraja, H.N. (2003).\  {\it Order Statistics}, Third  edition, John Wiley \& Sons, NY.

\item Durbin J., Stuart,  A. (1951).\ Inversions and rank correlation coefficients, {\it Journal of the Royal Statistical Society},  Ser. B {\bf 13} (2), 303--309.

\item Egorov, V. A., Nevzorov, V. B. (1984).\ Rate of convergence to the Normal law of sums of induced order statistics,
{\it Journal of Soviet Mathematics} (New York), {\bf 25}, 1139--1146.

\item Esscher, F. (1924).\  On a method of determining correlation from the ranks of the variates, {\it Skandinavisk Aktuarietidskrift}, {\bf 7}, 201--219.

{\it  Ann. Probab.}, {\bf 22}, 126--144.

\item Hauke, J.,  Kossowski, T. (2011).\ Comparison of values of Pearson's and Spearman's correlation coefficients on the same sets of data, {\it Quaestiones Geographicae}, {\bf 30} (2), 87--93.

\item Hoeffding, W. (1948).\ A class of statistics with asymptotically normal distribution, {\it Annals of Mathematical Statistics}, {\bf 19} (3),
293–-325.

\item Hotelling, H. (1953).\ New light on the correlatio ncoefficient and its transforms, Journal of the Royal Statistical Society: Series B (Statistical Methodology) {\bf 15} (2),   193–-232.

\item Fisher, R.A. (1921).\ On the ''probable error" of a coefficient of correlation deduced from a small sample, Metron, {\bf 1},  3--32.

\item Jensen, D.R. (1988).\  Semi-independence. In: Encyclopedia of Statistical Sciences, Volume 8, S. Kotz and N.L. Johnson (eds.),  358–-359. John Wiley and Sons, New York?

\item Kendall, M. G. (1970).\ {\it Rank Correlation Methods}, London, Griffin.

\item Omey, E., Gulk, S.Van (2008).\ Central limit theorems for variances and correlation coefficients, {\it Preprint}, doi: 10.13140/RG.2.2.21663.66727.

\item Shevlyakov, G.L., Vilchevski, N.O. (2002).\ {\it Robustness in Data Analysis: Criteria and Methods, Modern Probability and Statistics}, VSP, Utrecht.

\item Stepanov, A.V. (2025a).\ On Kendall's correlation coefficient, {\it Vestnik St. Petersburg University, Mathematics}, {\bf 58} (1), 71--78.

\item Stepanov, A. (2025b).\ Comparison of Correlation Coefficients,  {\it Sankhya A}, {\bf 87} (1),  191--218.

\item Walsh, J.E. (1969).\ Sample sizes for approximate independence between sample median and largest (or smallest) order statistics, {\it Austral.J. Statist.}, {\bf 11} (3), 120--122.

\item  Xu, W., Hou,  Y.,   Hung, Y. S., Zou, Y. (2013).\ A comparative analysis of Spearman's rho and Kendall's tau
in normal and contaminated normal models, {\it Signal Processing}, {\bf 93}, 261–-276.

}

\end{description}

\end{document}